\documentclass{amsart}
\usepackage{euler,amsmath,amssymb,amscd,amsfonts}
\usepackage{epsfig, graphicx,pinlabel}
\title[State polytope decomposition formula]{A state polytope decomposition formula}
\date{\today}
\author{Donghoon Hyeon and Jaekwang Kim}
\address[DH]{Department of Mathematics, POSTECH, Pohang, Gyungbuk 790-784, R. O. Korea}
\email{dhyeon@postech.ac.kr}
\address[JK]{Department of Mathematics, POSTECH, Pohang, Gyungbuk 790-784, R. O. Korea}
\email{kjk429@postech.ac.kr}
\newtheorem{theorem}{Theorem}[section]

\newtheorem{lemma}[theorem]{Lemma}
\newtheorem{proposition}[theorem]{Proposition}
\newtheorem{corollary}[theorem]{Corollary}
\theoremstyle{definition}
\newtheorem{defn}[theorem]{Definition}
\newtheorem{example}[theorem]{Example}

\def\st{\, | \,}

\def\mbb{\mathbb}

\def\a{\alpha}
\def\b{\beta}

\def\Sm{\Sigma}

\def\bR{\mbb R}
\def\bP{\mbb P}

\def\bN{\mbb N}
\def\bQ{\mbb Q}

\def\SL{{\mathrm SL}}
\def\GL{{\mathrm GL}}

\def\bP{\mathbb P}
\def\bG{\mathbb G}

\def\inj{\hookrightarrow}

\def\deg{\textup{deg} \, }

\def\aut{\textup{Aut}}

\def\inj{\hookrightarrow}

\def\Gr{Gr}

\def\bar{\overline}

\def\st{{\mathcal P}}

\input xy
\xyoption{all}

\input epsf
\def\bQ{\mathbb Q}

\def\M{\bar M}

\def\ii{\textup{\bf in}\,}
\newtheorem{remark}[theorem]{Remark}

\begin{document}

\begin{abstract} We give a decomposition formula for computing the state polytope of a reducible variety in terms of the state polytopes of its components: If a polarized projective variety $X$ is a chain of subvarieties $X_i$ satisfying some further conditions, then the state polytope of $X$ is the Minkowski sum of the state polytopes of $X_i$ translated by a vector $\tau$ which can be readily computed from the ideal of $X_i$. The decomposition is in the strongest sense in that the vertices of the state polytope of $X$ are precisely  the sum of vertices of the state polytopes of $X_i$ translated by $\tau$. We also give a similar decomposition formula for the Hilbert-Mumford index of the Hilbert points of $X$. We give a few examples of the state polytope and the Hilbert-Mumford index computation of reducible curves which are interesting in the context of the log minimal model program for the moduli space of stable curves.
\end{abstract}

\thanks{The first author was partially supported by the following grants funded by the government of Korea:
NRF grant 2011-0030044 (SRC-GAIA), NRF grant 2011-0005072 and the Korea Institute for
Advanced Study (KIAS) grant. The second author was partially supported by the following grants funded by the government of Korea:
NRF grant 2011-0030044 (SRC-GAIA), NRF grant 2011-0005072. }
\maketitle

\section{Introduction}
The state polytope of an ideal encodes much information about the scheme it defines. Let $V$ be a vector space over an algebraically closed field $k$ of characteristic zero. Given a rational representation $W$ of $\SL(V)$ and a maximal torus $T \subset \SL(V)$, Kempf \cite[\S3]{Kempf} defined the {\it state} of $w\in W$ (with respect to $T$) to be the set of the characters $\chi \in X(T)$ such that $w_\chi \ne 0$ where $w_\chi$ is the projection of $w$ in the weight space $W_\chi$.  Given a projective variety $X \subset \bP(V)$ and a choice of homogeneous coordinates, Bayer and Morrison in \cite{BM} defined the {\it $m$th state polytope} of $X \subset \bP(V)$ to be the convex hull of the states of (any affine point over) the $m$th Hilbert point $[X]_m \in \bP\left(\bigwedge^{Q(m)}S^mV^*\right)$ defined:
\[
\left[(I_X)_m \to S^mV^*  \right] \in \Gr_{Q(m)}S^mV^*) \inj  \bP\left(\bigwedge^{Q(m)}S^mV^*\right)
\]
 where $I_X$ is the saturated homogeneous ideal of $X$, $Q(t)$ is the dimension of the $t$th graded piece $(I_X)_t$,  and $m$ is a positive integer bigger than or equal to the Castelnuovo-Mumford regularity of $X$. The state polytope of a homogeneous ideal is the state polytope of the projective variety it defines. The relation between the state polytopes and
the Gr\"obner theory is described by:

 \begin{theorem}\cite[Theorem~3.1]{BM}\label{T:BM} \, There is a natural one-to-one correspondence between the initial ideals and the vertices of the state polytope.
 \end{theorem}

 Let $T$ be the maximal torus of $\SL(V^*)$ diagonalized by $x_0, \dots, x_n$. Then by considering the $T$-weight space decomposition of $\bigwedge^{Q(m)}S^mV^*$, one can naturally associate the characters in the state of a Hilbert point $\left[(I_X)_m \to S^mV^*\right]$ and the {\it monomials} $x^{\a(1)}\wedge \cdots \wedge x^{\a(Q(m))}$ whose associated Pl\"ucker coordinate does not vanish at $[X]_m$.  See \cite{MS} for this correspondence (Section~3.1) as well as for a very nice exposition on Kempf's theory of the worst one-parameter subgroup and basics on state polytopes. In particular, the trivial character corresponds to the {\it barycenter} of which coordinates are all $\frac{mQ(m)}{\dim V}$.

In view of this correspondence, we take the following definition of the state polytope \cite[Formula~(2.7)]{Sturmfels}:
\begin{defn}
Given a homogeneous ideal $I \subset S = k[x_0,\dots, x_n]$ and $m \ge \textup{reg}(I)$, the {\it $m$th state polytope} is  defined and denoted by
\begin{equation}\label{E:st}
\st_m(I) := \textup{Conv}\left\{\left. \sum_{x^{\a} \in \ii_{\prec}(I)_m} \a  \, \right|  \mbox{$\prec$ is a monomial order}  \right\}.
\end{equation}
\end{defn}
Here, Conv means taking the convex hull and $\textup{reg}(I)$ is the Castelnuovo-Mumford regularity of $I$, and the condition  $m \ge \textup{reg}(I)$ ensures that the $m$th Hilbert point of the scheme $X$ cut out by $I$ is well defined.

Apparent from the definition is that the state polytope can be computed from the universal Gr\"obner basis. More importantly,  $\st_m(I_X)$ determines the semistability of the Hilbert point $[X]_m$  with resect to the chosen basis. This is a direct consequence of the  Mumford's numerical criterion \cite{GIT}: $[X]_m$ is $T$-semistable (resp. $T$-stable) if and only if the state polytope in $X(T)$ (resp. interior of the polytope) contains the trivial character. This condition is equivalent to the state polytope (\ref{E:st}) (resp. its interior) containing the barycenter, where the coordinates $x_0,\dots, x_n$ are chosen so that $T$ acts on them via characters i.e. they diagonalize the $T$ action. The upshot is that, by computing the universal Gr\"obner basis (with a computer algebra system if and when convenient), one can determine the semistability with respect to the given coordinates i.e. with respect to the associated maximal torus $T$. If $[X]_m$ is $T$-unstable for some $T$, then $[X]_m$ is GIT unstable.

Of course, to prove the semistability of $[X]_m$, one has to prove its $T$-semistability for {\it all} maximal torus $T$, so in that regard the state polytope formulation of GIT semistability may not seem too much of a help. But in the special case when $X \subset \bP(V)$ and $V$ is a multiplicity-free representation of a linearly reductive subgroup $\Gamma$ of $\aut(X)$, the (semi)stability of $[X]_m$ is equivalent to the $T$-(semi)stability with respect to any maximal torus $T$ that preserves the $\Gamma$-irreducibles of $V$ \cite[Proposition~4.7]{MS}. This is the key idea of Morrison and Swinarski that allowed them to prove the $m$-Hilbert semistability of various curves for small $m$ (\S7, ibid).  It was also the starting point for Alper, Fedorchuk and Smyth to obtain their results on the $m$-Hilbert semistability of canonical and bicanonical images of generic smooth curves for $m\ge 2$ \cite{AFS}, which should prove essential in carrying out the log minimal model program for $\overline M_g$ ({\it the Hassett-Keel program}). In fact, Alper, Fedorchuk and Smyth do not rely on the state polytope technique: Instead, they work out {\it by hand} a collection of basis members for the (bi)canonical system and deduce from them the semistability directly, which is in every manner very impressive.

\

Inspired from the exciting developments, we shall consider in this article how one can more efficiently compute the state polytope of certain reducible varieties. More precisely, we give a formula for the state polytope of a variety in terms of the state polytopes of its subvarieties. We say that $X$ is a {\it chain of subvarieties $X_i$} if  $X = \cup_{i=1}^\ell X_i$ and $X_i$ meets $X_j$ when and only when $|i-j|=1$.

\begin{theorem}\label{T:main} Let $X$ be a chain of subvarieties $X_1, \dots, X_\ell$ defined by a homogeneous ideal $I_X = \cap_i I_{X_i}$. Suppose that there is a homogeneous coordinate system $x_0,\dots,x_n$ and a sequence $n_{0}=0 < n_1 < \cdots < n_\ell = n$ such that
\[
X_i \subset \{x_0 = \cdots = x_{n_{i-1}-1} = 0, x_{n_i+1} = x_{n_i+2} = \cdots = x_n = 0\}.
\]
Then the state polytope of $X$ is given by the following decomposition formula
\small\begin{equation}\label{E:main}
\st_m(I_X) = \sum_{i=1}^\ell \st_m( I_{X_i}\cap k[x_{n_{i-1}},\cdots,x_{n_i}]) + \sum_{i=1}^{\ell-1} \st_m(T_i\cap k[x_{n_{i-1}}, \dots, x_n])  \tag{$\dagger$}
\end{equation} \normalsize
where  $T_i = \langle x_{n_{i-2}},\dots,x_{n_{i-1}-1}\rangle \langle x_{n_i+1},\dots,x_n \rangle$ for $2 \le i \le \ell-1$, and $T_1 = \langle x_{n_1+1},x_{n_1+2},\dots,x_n\rangle$ and $T_{\ell} = \langle x_{n_{\ell-2}},x_{n_{\ell-2}+1},\dots,x_{n_{\ell-1}-1}\rangle$.
\end{theorem}

\begin{remark}\label{rem:main} \begin{enumerate}
\item Here, $\st_m( I_{X_i}\cap k[x_{n_{i-1}},\cdots,x_{n_i}])$ is regarded as a convex polytope in the subspace
\[
\{(a_0,\dots,a_n) \in \bR^{n+1} \, | \, a_0 = \cdots = a_{n_{i-1}-1} = 0, a_{n_i+1} = a_{n_i+2} = \cdots = a_n = 0\}.
\]
Similarly, $\st_m(T_i\cap k[x_{n_{i-1}}, \dots, x_n])$ is also regarded as a convex polytope in the relevant vector subspace.

\item
Note that the second term of (\ref{E:main}) is zero dimensional since $T_i$ are monomial ideals. We shall reserve the letter $\tau$ to denote it.
\end{enumerate}
\end{remark}

In fact, the polytope decomposition in Theorem~\ref{T:main} is sharp in the following sense:

\begin{corollary}\label{C:vertex} Retain notations from Theorem~\ref{T:main}. Let $\mathcal V_i$ denote the set of vertices of $\st_m(I_{X_i}\cap k[x_{n_{i-1}},\dots,x_{n_i}])$, $i = 1, \dots, \ell$. Then the vertices of $\st_m(I_X)$ are precisely
\[
\left\{ \left. \tau + \sum_{i=1}^\ell v_i \, \right| \, v_i \in \mathcal V_i \right\}.
\]
\end{corollary}

The Hilbert-Mumford index can also be computed by a similar decomposition formula. To be consistent with our main reference \cite{HHL}, we state the formula in terms of the {\it dual Hilbert point}
\[
[X]_m^\star = [S^mV^*\to S^mV^*/I_m] \in \Gr^{P(m)}S^mV^* \inj \bP\left(\bigwedge^{P(m)}S^mV^*\right).
\]
where $P(t) \in \bQ[t]$ is the Hilbert polynomial of $X$.

\begin{proposition}\label{P:mu-reduction} Let $X$ be as in Theorem~\ref{T:main} and $\rho:\mathbb{G}_m\to \GL_{n+1}$ be a 1-parameter subgroup of $\GL_{n+1}$ diagonalized by $\{x_0,\dots,x_n\}$ with weights $(r_0,\cdots,r_n)$ and $\rho_i$ be the restriction of $\rho$ to $\GL(kx_{n_{i-1}}+\cdots+kx_{n_i})$.  Then the Hilbert-Mumford index $\mu([X]_m^\star,\rho)$ of the $m$th Hilbert point of $X$ with respect to $\rho$ is given by
\small
$$
\mu([X]_m^\star,\rho)=\sum_{i=1}^{\ell}\mu([X_i]_m^\star,\rho_i)-\sum_{i=1}^{\ell}\bigg(\frac{mP_i(m)}{n_{i}-n_{i-1}+1}\sum_{k=n_{i-1}}^{n_i}r_k\bigg)
+ \frac{m P(m)}{n+1}\sum_{i=0}^n r_i +m\sum_{i=1}^{\ell-1}r_{n_i}
$$
\normalsize
where $P(m)$ is the Hilbert polynomial of $I_X \subset k[x_0,\dots,x_n]$ and  $P_i(m)$, the Hilbert polynomial of $I_{X_i} \cap k[x_{n_{i-1}}, \dots, x_{n_i}]$ regarded as an ideal in $k[x_{n_{i-1}},\dots,x_{n_i}]$.
\end{proposition}

We shall give  proofs of Theorem~\ref{T:main} and Proposition~\ref{P:mu-reduction} in \S\ref{S:st-reduction}. Corollary~\ref{C:vertex} will be established in \S\ref{S:reduction-in}.
\section{Basic examples}
Before proving the main results, we shall give a few basic examples at the far ends of the spectrum, namely, monomial ideals and hypersurfaces (plane curves, to be more specific).
\begin{example}[Monomial ideals] Let $X = X_1\cup X_2 \cup X_3 \subset \bP^3$ be a chain of $\bP^1$'s:
\[
X_1 = \{x_2 = x_3 = 0\}, \quad X_2 = \{x_0 = x_3 = 0\}, \quad X_3=\{x_0=x_1=0\}.
\]
Then $I_{X_1} \cap k[x_0,x_1] = 0$  and it does not contribute to the state polytope. The other components do not contribute for the same reason, but the mixed terms
 \[
 \{x_0x_2, x_0x_3\} \cup \{x_0x_3,x_1x_3\}
 \]
are precisely the monomial generators of the ideal of $X$, and the formula (\ref{E:main}) holds.

\

 More generally, if $X = \cup X_i$ as in Theorem~\ref{T:main} and each $X_i$ are defined by monomials $M_{i\a}$,
\[
\st_m(I_{X_i} \cap k[x_{n_{i-1}},\dots, x_{n_i}]) = \left\{\sum_\a \log_x M_{i\a}\right\}
\]
where $\log_x x^\a = \a$. The formula (\ref{E:main}) implies that
\[
\st_m(I_X) = \left\{\sum_i \sum_\a \log_x M_{i\a} + \sum_{x^\b \in \sum_i T_i} \b\right\}.
\]
This is the same as
\[
\st_m(I_X) = \left\{\sum_{x^\a \in (I_X)_m} \a \right\}
\]
since
\[
I_X = \cap_{i} I_{X_i} = \cap_i \langle \{M_{i\a}\}_\a, x_0,\dots,x_{n_{i-1}-1}, x_{n_i+1}, \dots, x_n\rangle = \langle \{M_{i\a}\}_{i,\a} \rangle + \sum_i T_i.
\]
Note that the inclusion $M_{i\a} \in \cap_i I_{X_i}$ follows from the assumption that $M_{i\a}$ is not a power of $x_{n_i}$ or of $x_{n_{i-1}}$.
\end{example}

\begin{example}[Plane curves] We consider a simple example of two plane curves
\[
E_1 = \{b^2c=a(a-c)(a-2c), d=e=0\} \quad \mbox{and} \quad E_2 = \{d^2c=e^2(e+c),a=b=0\}
\]
meeting in one node. Let $C$ denote the union of $E_1$ and $E_2$. The $3$rd state polytope of $E_1$ has three vertices
\begin{equation}\label{E:E1}
\{(3, 0, 0, 0, 0), (1, 0, 2, 0, 0), (0, 2, 1, 0, 0)\}
\end{equation}
and that of $E_2$ has two vertices
\begin{equation}\label{E:E2}
\{(0, 0, 1, 2, 0), (0, 0, 0, 0, 3)\}.
\end{equation}
Indeed, since $E_i$ are hypersurfaces (in suitable linear subspaces) of degree three, their 3rd state polytopes are precisely the Newton polytopes.

To compute the state polytope of $C$, we first take the sums of a point from $(\ref{E:E1})$ and a point from $(\ref{E:E2})$:
$$\{(3,0,1,2,0),
(3,0,0,0,3),
(1,0,3,2,0),
(1,0,2,0,3),
(0,2,2,2,0),
(0,2,1,0,3)\}.$$
By the decomposition formula, the vertices themselves are these translated by $\tau$ which is the sum of the exponent vectors of
\[
T = \{a^2d, abd,acd,ad^2,ade,a^2e,abe,ace,ae^2,b^2d,bcd,bd^2,bde,b^2e,bce,be^2\}.
\]
We compute $\tau = (11,11,4,11,11)$, and hence the vertices of $\st_3(I_C)$ are:
\[\begin{array}{c}
(14,11,5,13,11),
(14,11,4,11,14),
(12,11,7,13,11),\\
(12,11,6,11,14),
(11,13,6,13,11),
(11,13,5,11,14).
\end{array}
\]

This agrees with the direct computation of the $3$rd state polytope of the ideal of $C$:
\[
\langle b e,a e,b d,a d,-c d^{2}+e^{3}+e^{2},a^{3}-3 a^{2} c-b^{2} c+2 a
     c^{2} \rangle
 \]
 Here we demonstrate the output of the Macaulay 2 \cite{M2} package \verb"StatePolytope" written by D. Swinarski.

\footnotesize
\begin{verbatim}
i2 : R = QQ[a,b,c,d,e];
i3 : I1 = ideal(b^2*c - a*(a-c)*(a-2*c))
              3     2     2        2
o3 = ideal(- a  + 3a c + b c - 2a*c )
o3 : Ideal of R
i4 : I2 = ideal(d^2*c-e^2*(e+1))
              2    3    2
o4 = ideal(c*d  - e  - e )
o4 : Ideal of R
i5 : I = intersect(I1+ideal(d,e),I2+ideal(a,b))
                                     2    3    2   3     2     2        2
o5 = ideal (b*e, a*e, b*d, a*d, - c*d  + e  + e , a  - 3a c - b c + 2a*c )
o5 : Ideal of R
i6 : statePolytope(3,I)
LP algorithm being used: "cddgmp".
polymake: used package cddlib
  Implementation of the double description method of Motzkin et al.
  Copyright by Komei Fukuda.
  http://www.ifor.math.ethz.ch/~fukuda/cdd_home/cdd.html

VERTICES
1 14 11 5 13 11
1 14 11 4 11 14
1 12 11 6 11 14
1 11 13 5 11 14
1 12 11 7 13 11
1 11 13 6 13 11


o6 = {{14, 11, 5, 13, 11}, {14, 11, 4, 11, 14}, {12, 11, 6, 11, 14},
     ---------------------------------------------------------------
      {11, 13, 5, 11, 14}, {12, 11, 7, 13, 11}, {11, 13, 6, 13, 11}}
\end{verbatim}\normalsize

\end{example}

\begin{example}\label{Ex:rhamphoid} This example is non-trivial compared to the previous two. We shall consider a particular genus four curve with a genus two tail. Let $R$ be a rational curve with a rhamphoid cusp. It is of arithmetic genus two and admits a $\bG_m$ action with two fixed points one of which is the cusp. The action comes from the automorphism of its normalization. Let $C$ be a genus two curve obtained by attaching two copies $R_1, R_2$  of $R$ at the smooth fixed points, say $p_i \in R_i$, $i=1,2$.
\begin{figure}[!htb]\labellist \small\hair 2pt
\pinlabel $R_1$ at 85 625
\pinlabel $y^2=x^5$ at 136 685
\pinlabel $p$ at 215 635
\pinlabel $R_2$ at 350 625
\endlabellist
  \begin{center}
    \includegraphics[width=2.5in]{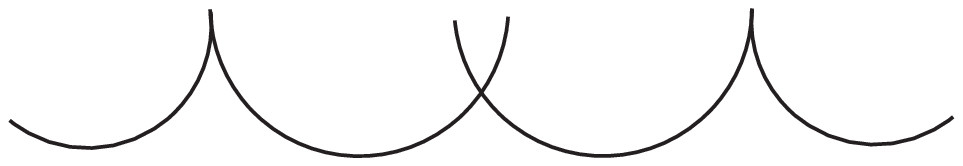}
  \end{center}
   \caption{} \label{F:ramphoid-cusp}
\end{figure}
We bicanonically embed $C$ in $\bP^8$ and consider its state polytope.
Note that $R_1$ can be parametrized by
\[
\begin{array}{lll}
\bP^1 & \rightarrow & \bP^8\\
\left[s,t\right] & \mapsto & \left[s^6, s^4t^2, s^2t^4, st^5, t^6, {\bf 0}_4\right]
\end{array}
 \]
which has a rhamphoid cusp at $[1,0, \dots, 0]$ and $p:=[{\bf 0}_4, 1, {\bf 0}_4]$ is fixed under the automorphism where ${\bf 0}_i = (\underbrace{0, \dots, 0}_i)$. $R_2$ is parametrized similarly, from which we can compute its defining saturated ideal. The $6$th state polytope of $I_{R_1}\cap k[x_0,\dots,x_4]$ has $51$ vertices
\small
\[
 \{(216,191,206,206,231,{\bf 0}_{4}),(216,191,191,236,216,{\bf 0}_{4}),\cdots,(181,248,210,180,231,{\bf 0}_{4})\}
\]\normalsize
which are obtained by Swinarski's \verb"StatePolytope".

The state polytope of $I_{R_2}\cap k[x_4,\dots, x_8]$ is the mirror flip of the above with respect to $x_4$. By the decomposition formula, the $6$th state polytope $P$ of $C$ is $\tau$-translate of the convex hull of $51\cdot 51 = 2601$ vertices obtained by choosing one from each state polytope and adding them up. Being the state polytope of a product of ideals of linear subspaces, $\tau$ can be readily computed either by hand or a computer algebra system:
\[
\tau  =(1750,1750,1750,1750,1504,1750,1750,1750,1750).
 \]
The barycenter of this state polytope is $v=(1956,1956,\dots,1956)$. How can we utilize the decomposition formula to check whether $v$ is in the state polytope? A natural thing to do is to try to decompose $v-\tau$ into the two symmetric points $v_1=(206,206,206,206,226,{\bf 0}_4)$ and $v_2=({\bf 0}_4,226,206,206,206,206)$, and see if $v_i$ is contained in $6$th state polytope of $R_i$.  But this can be easily checked by using the function \verb"contains" of the polyhedra package in Macaulay 2 as follows.
\begin{verbatim}
i3 : R=QQ[a,b,c,d,e];
i4 : Q=QQ[s,t];
i5 : f=map(Q,R,{s^6,s^4*t^2,s^2*t^4,s*t^5,t^6});
o5 : RingMap Q <--- R
i6 : I=ker f;
o6 : Ideal of R
i7 : L=statePolytope(6,I);
i8 : P1=convexHull transpose matrix L;
i9 : v1=transpose matrix {{206,206,206,206,226}};
              5        1
o9 : Matrix ZZ  <--- ZZ
i10 : contains(P1,v1)
o10 = true
\end{verbatim}\normalsize
Hence, we conclude that $v$ is contained in $P$.
In fact, the ideal of $C$ is simple enough so that its state polytope can be directly computed by using the state polytope package of Macaulay 2, which agrees with the result obtained above.
\end{example}

The assumption on the existence of the coordinate system as in  Theorem~\ref{T:main} may seem quite restrictive, but we shall see that it is satisfied by an important class of varieties namely, the pluricanonical images of the generic members of the boundary of $\M_g$. We shall give a few interesting examples in this vein in \S\ref{S:applications}.


\section{Decomposition formula for initial ideals}\label{S:reduction-in}
First we shall prove a key lemma on initial ideals from which the main theorem follows with some simple observations regarding the monomial orders.
Let $Y$ and $Z$ be closed subvarieties in $\mathbb{P}^n$ defined by homogeneous ideals $I_Y$ and $I_Z$ of $k[x_0,\dots,x_n]$ respectively, and let $X$ be the projective variety defined by $I_X := I_Y\cap I_Z$. Suppose that with respect to the homogeneous coordinate system $x_0,x_1,\cdots,x_n$, the subvarieties $Y$ and $Z$ are contained in linear subspaces as follows:
\[
\begin{array}{l}
Y\subset L_1 := \{x_0=\cdots=x_{l-1}=0\}\cong\mathbb{P}^{n-l} \\
Z\subset L_2 := \{x_{l+1}=\cdots=x_n=0\}\cong\mathbb{P}^l.
\end{array}
\]
In particular, $Y\cap Z = \{ p \}$
where $p$ is the unique point in $L_1\cap L_2$ whose coordinates are all zero except $x_l$.

\begin{lemma}\label{L:in} Let $\prec$ be a monomial order. The initial ideal of $I$ with respect to $\prec$ is given by
\begin{equation}\label{E:in}
\ii_\prec(I_X) = \langle \ii_\prec(I_Y\cap k[x_l,\cdots,x_n])\rangle + \langle\ii_\prec(I_Z\cap k[x_0,\cdots,x_l]) \rangle+ T \tag{$\dagger\dagger$}
\end{equation} \normalsize
where $T = \langle x_0,\dots,x_{l-1} \rangle \langle x_{l+1},\dots,x_n\rangle$. Note that $\ii_\prec(I_Y\cap k[x_l,\cdots,x_n])$  is computed as an ideal of $k[x_l, \cdots, x_n]$ and $\langle \ii_\prec(I_Y\cap k[x_l,\cdots,x_n])\rangle$ is the ideal in $k[x_0,\cdots,x_n]$ it generates. A similar statement for $Z$ is noted.
\end{lemma}

\begin{proof}
  We first note that $T$ is contained in $\ii_\prec(I_X)$ since it is a monomial ideal contained in $I_X$.

  Let $x^\alpha=\prod x_i^{\alpha_i}$ be a  monomial in   $\ii_\prec(I_X)_m$ i.e. $x^\alpha=\ii_\prec(f)$ for some $f\in I_X$. If $x^\alpha\notin T$, then $x^\alpha$ is contained in $k[x_0,\cdots,x_l]$ or $k[x_l,\cdots,x_n]$. If $x^\alpha \in k[x_0,\cdots,x_l]$, then $\ii_\prec(g)=x^\alpha$ where $g(x_0,\dots,x_l) = f(x_0,\cdots,x_l,0,\cdots,0)$.  But $g\in I_Z\cap k[x_0,\cdots,x_l]$, so $x^\alpha$ is contained $\ii_\prec(I_Z\cap k[x_0,\cdots,x_l])_m$. Similarly, if $x^\alpha\in k[x_l,\cdots,x_n]$, then $x^\alpha\in \ii_\prec(I_Y\cap k[x_l,\cdots,x_n])_m$. This proves that the left hand side is contained in the right.

To see the other inclusion,  suppose that $x^\alpha\in \ii_\prec(I_Y\cap k[x_l,\cdots,x_n])_m$ i.e.,  $x^\alpha=\ii_\prec(f)$ for some $f\in I_Y\cap k[x_l,\cdots,x_n]$. Monomials of $k[x_l,\dots,x_n]$ that do not vanish at $p$ are of the form $x_l^i$. But since $f$ is homogeneous and vanishes at $p \in Y$, it cannot have the term $x_l^{\deg(f)}$, so each term of $f$ is divisible by $x_i$ for some $i>l$. This implies that $f$ vanishes on $Z$ so that $f\in I_X$ and $x^\alpha\in \ii_\prec(I_X)$. Hence $\ii_\prec (I_Y\cap k[x_l,\cdots,x_n]) \subset \ii_\prec(I_X)$, and by a similar argument $\ii_\prec(I_Z\cap k[x_0,\cdots,x_l]) \subset \ii_\prec(I_X)$.
\end{proof}

We obtain the following corollary by induction.

\begin{corollary}\label{C:in} Retain the notations from the main Theorem~\ref{T:main}. For any monomial order $\prec$ on $k[x_0,\dots,x_n]$, we have
\[
\begin{array}{cll}
\ii_\prec(I_X) & = & \sum_{i=1}^{\ell}\langle \ii(I_{X_i}\cap k[x_{n_{i-1}},\cdots,x_{n_i}]) \rangle + \sum_{i=1}^{\ell-1}T_i \\
& = & \sum_{i=1}^{\ell}\langle \ii(I_{X}\cap k[x_{n_{i-1}},\cdots,x_{n_i}]) \rangle + \sum_{i=1}^{\ell-1}T_i.
\end{array}
\]
\end{corollary}

\begin{proof} Suppose that the formula holds for $\ell-1$. Regarding $X$ as the union of $Z:=\cup_{i=1}^{\ell-1}X_i$ and $Y := X_\ell$ and applying Lemma~\ref{L:in} and the induction hypothesis, we obtain
\[
\begin{array}{lll}
\ii_\prec(I_X)
& = &  \sum_{i=1}^{\ell-1}\langle \ii(I_{X_i}\cap k[x_{n_{i-1}},\cdots,x_{n_i}]) \rangle + \sum_{i=1}^{\ell-2}T_i \\
& & + \langle \ii_\prec(I_Z\cap k[x_{n_{\ell-1}}, \dots, x_n])\rangle  + T\\
\end{array}
\]
where $T = \langle x_0,\dots, x_{n_{\ell-1}-1}\rangle\langle x_{n_{\ell-1}+1}, \dots, x_n\rangle$ may be replaced by $$T_{\ell-1}=\langle x_{n_{\ell-2}},\dots, x_{n_{\ell-1}-1}\rangle\langle x_{n_{\ell-1}+1}, \dots, x_n\rangle$$ since
$\langle x_0,\dots, x_{n_{\ell-2}-1}\rangle\langle x_{n_{\ell-1}+1}, \dots, x_n\rangle$ is contained in $\sum_{i=1}^{\ell-2} T_i$.
\end{proof}

Now we can give
\begin{proof}[Proof of Corollary~\ref{C:vertex}]
This is an immediate consequence of Theorem~\ref{T:BM} and Corollary~\ref{C:in}.
\end{proof}

\begin{lemma}\label{L:standard}(Standard monomials)
Let $\Sigma_{X,m}$, $\Sigma_{Y,m}$ and $\Sigma_{Z,m}$ be defined
\[
\begin{array}{lll}
\Sigma_{X,m} & =  & \{\mbox{monomials in $k[x_0, \dots, x_n]_m$}\} \setminus \ii_\prec(I_X)_m  \\
\Sigma_{Y,m}  & =  & \{\mbox{monomials in $k[x_l, \dots, x_n]_m$}\} \setminus \ii_\prec(I_Y\cap k[x_l,\cdots,x_n])_m  \\
\Sigma_{Z,m} & = & \{\mbox{monomials in $k[x_0,\cdots,x_l]_m$}\} \setminus \ii_\prec(I_Z\cap k[x_0,\cdots,x_l])_m.
\end{array}
\]
Then $\Sigma_{X,m}=\Sigma_{Y,m}\cup \Sigma_{Z,m}$ and $\Sigma_{Y,m}\cap \Sigma_{Z,m}=\{x_l^m\}$.
\end{lemma}

\begin{proof}
Let $x^\a$ be a degree $m$ monomial and let $T$ be as in Lemma~\ref{L:in}. Since $T \subset I_X$, if $x^\a$ is not in $\ii_\prec(I_X)$, it is necessarily in $k[x_0,\dots,x_l]$ or in $k[x_l,\dots,x_n]$. Such a monomial $x^\a$ is in the initial ideal if and only if it is in $\ii_\prec(I_Y\cap k[x_l,\cdots,x_n])_m$ or in $\ii_\prec(I_Z\cap k[x_0,\cdots,x_l])_m$, and the equality $\Sigma_{X,m} = \Sigma_{Y,m}\cup \Sigma_{Z,m}$ follows.
That $\Sigma_{Y,m}\cap \Sigma_{Z,m}=\{x_l^m\}$ is clear from the fact that $x_l^m$ is the only degree $m$ monomial not vanishing at $p$, and thus not contained in any of the initial ideals involved in the discussion.
\end{proof}

\section{State polytope decomposition formula}\label{S:st-reduction}
We prove our main results Theorem~\ref{T:main} and Proposition~\ref{P:mu-reduction} in this section.

\begin{proof}[Proof of Theorem~\ref{T:main}] Surely, it suffices to prove Theorem~\ref{T:main} for the $\ell = 2$ case, as the general case would follow from it by a simple induction as in the proof of Corollary~\ref{C:in}. So, let $X = Y\bigcup Z$ and $T$ be as in Lemma~\ref{L:in}. We shall prove that
\begin{equation}\label{E:st-reduction}
\st_m(I_X)=\st_m(I_Y\cap k[x_l,\cdots,x_n])+\st_m(I_Z\cap k[x_0,\cdots,x_l])+\tau
\end{equation}
where $\tau=\underset{x^\alpha\in T}\sum\alpha$.

Let $\sum_{x^\alpha\in \ii_\prec(I_X)_m}\alpha$ be a vertex of $\st_m(I_X)$ induced by some monomial order $\prec$ on $k[x_0,\cdots,x_n]$.
Define
 $$\Sigma_{Y,m}^c=\ii_\prec(I_Y\cap k[x_l,\cdots,x_n])_m, \quad \Sigma_{Z,m}^c=\ii_\prec(I_Z\cap k[x_0,\cdots,x_l])_m. $$
By Lemma~\ref{L:in}, we have
\[
\ii_\prec(I_X) =\langle \ii_\prec(I_Y\cap k[x_l,\cdots,x_n]) \rangle +  \langle \ii_\prec(I_Z\cap k[x_0,\cdots,x_l]) \rangle + T
\] which implies
\[
\sum_{x^\alpha\in \ii_\prec(I_X)_m}\alpha=\sum_{x^\alpha\in \Sigma_{Y,m}^c}\alpha+\sum_{x^\alpha\in \Sigma_{Z,m}^c}\alpha+\sum_{x^\alpha\in T}\alpha.
\]
Since $\sum_{x^\alpha\in \Sigma_{Y,m}^c}\alpha$ and $\sum_{x^\alpha\in \Sigma_{Z,m}^c}\alpha$ are vertices of $\st_m(I_Y\cap k[x_l,\cdots,x_n])$ and $\st_m(I_Z\cap k[x_0,\cdots,x_l])$ respectively, $\sum_{x^\alpha\in \ii_\prec(I_X)_m}\alpha$ is contained in the right hand side of (\ref{E:st-reduction}).

Conversely, let $\alpha_1$ be a vertex of $\st_m(I_Y\cap k[x_l,\cdots,x_n])$ induced by the monomial order $\prec_1$ on $k[x_l,\cdots,x_n]$ and $\alpha_2$, a vertex of $\st_m(I_Z\cap k[x_0,\cdots,x_l])$ induced by the monomial order $\prec_2$ on $k[x_0,\cdots,x_l]$. We claim that there is a monomial order $\prec$ on $k[x_0,\cdots,x_n]$ that induces the initial ideals with respect to the given orders $\prec_1$ on $k[x_l,\cdots,x_n]$ and $\prec_2$ on $k[x_0,\cdots,x_l]$.
There are vectors $v \in \bN^{n-l+1}$ and $ v' \in \bN^{l+1}$  such that (cf. \cite[Proposition~1.11]{Sturmfels}).
\[
\begin{array}{lll}
\ii_{\prec_v}(I_Y\cap k[x_l,\cdots,x_n]) &= &\ii_{\prec_1}(I_Y\cap k[x_l,\cdots,x_n]) \\
\ii_{\prec_{v'}}(I_Z\cap k[x_0,\cdots,x_l])& = &\ii_{\prec_2}(I_Z\cap k[x_0,\cdots,x_l])
\end{array}
\]
where $\prec_v$ and $\prec_{v'}$ are the weight orders given by $v$ and $v'$, respectively. In general, an integral vector only gives rise to a partial order, but the choice of $v$ and $v'$ is made such that they give total orders. This is possible, again by \cite[Proposition~1.11]{Sturmfels}.
By modifying the first entry of $v$ (without affecting the order it defines), we may assume that it equals the last entry of $v'$. For instance, we may simply add $v'_l - v_1$ to all entries of $v$, which does not change the monomial order. Let $w = (v'_1, \dots, v'_l=v_1, v_2, \dots, v_{n-l+1})$. In general, this only defines a partial order in which case we may employ any tie-breaking device, for instance the Lex order, to define a total order: declare $M = \prod_{i=0}^n x_i^{\alpha_i} \prec M' = \prod_{i=0}^n x_i^{\beta_i}$ if
\begin{enumerate}
\item[(i)] $w. \alpha < w. \beta$; or
\item[(ii)] $w. \alpha = w. \beta$ and $M \prec' M'$
\end{enumerate}
where $\prec'$ is the chosen tie-breaking. Note that $\ii_\prec(I_Y\cap k[x_l,\dots,x_n]) = \ii_{\prec_1}(I_Y\cap k[x_l,\dots,x_n])$ since $\prec$ induces the weight order given by $v$ on $k[x_l,\dots,x_n]$. Similar statement holds for $I_Z$ and $\prec_2$.
Then by Lemma~\ref{L:in}, we have
\[
\ii_\prec(I_X) = \langle\ii_{\prec_1}(I_Y\cap k[x_l,\cdots,x_n])\rangle + \langle\ii_{\prec_2}(I_Z\cap k[x_0,\cdots,x_l])\rangle + T
\]
and this shows that $\alpha_1+\alpha_2+\tau=\sum_{x^\alpha\in \ii_{\prec}(I_X)_m}\alpha\in \st_m(I_X)$.
This completes the proof of Theorem~\ref{T:main}.
\end{proof}

We now prove the decomposition formula for the Hilbert-Mumford index.

\begin{lemma}\label{L:mu-reduction}
Let $X, Y, Z$ be as in Lemma~\ref{L:in}, $\rho:\mathbb{G}_m\to \GL_{n+1}$ be a 1-parameter subgroup diagonazlied by $\{x_0, \dots, x_n\}$ with weights $(r_0,r_1,\cdots,r_n)$, and $\rho',\rho''$ be the restrictions of $\rho$ to $\GL(kx_l+\cdots+kx_n)$ and $\GL(kx_0+\cdots+kx_l)$, respectively. Then the Hilbert-Mumford index $\mu([X]_m^\star,\rho)$ of the $m$th Hilbert point of $X$ with respect to $\rho$ is given by
$$
\begin{array}{cll}
\mu([X]_m^\star,\rho) & = &  \mu([Y]_m^\star,\rho')+\mu([Z]_m^\star,\rho'') -\frac{mP_Z(m)}{n+1-l}\sum_{i=l}^{n}r_i-\frac{mP_Y(m)}{l+1}\sum_{i=0}^{l}r_i\\
& & + \frac{mP_X(m)}{n+1}\sum_{i=0}^nr_i + mr_l
\end{array}
$$
where $$P_Y(m)=\dim_k (k[x_l,\cdots,x_n]/\ii_{\prec_{\rho'}}(I_Y\cap k[x_l,\cdots,x_n]))_m$$ and $$P_Z(m)=\dim_k(k[x_0,\cdots,x_l]/\ii_{\prec_{\rho''}}(I_Z\cap k[x_0,\cdots,x_l]))_m$$ are the Hilbert polynomials of $Y$ and $Z$ regarded as closed subvarieties of  \linebreak
$\{x_0 = \cdots = x_{l-1} = 0\} \simeq \bP^{n-l}$ and $\{x_{l+1} = \cdots = x_n = 0\} \simeq \bP^l$, respectively.
\end{lemma}

\begin{proof}
Let $\Sigma_{X,m}$, $\Sigma_{Y,m}$ and $\Sigma_{Z,m}$ be as in Lemma~\ref{L:standard}. By \cite{HHL}, the Hilbert-Mumford index can be computed by the formula
\begin{equation}\label{E:HHL}
\mu([X]_m^\star,\rho)=-\underset{x^\alpha\notin \ii_{\prec_\rho}(I_X)_m}\sum wt(x^\alpha) + \frac{mP_X(m)}{n+1}\sum_{i=0}^n r_i.
\end{equation}
Since $\Sigma_{X,m}=\Sm_{Y,m}\cup \Sm_{Z,m}$ and $\Sm_{Y,m}\cap \Sm_{Z,m}=\{x_l^m\}$ by Lemma~\ref{L:standard}, we have
$$
-\underset{x^\alpha\notin in(I_X)_m}\sum wt(x^\alpha)=-\underset{x^\alpha\in \Sm_{Y,m}}\sum wt(x^\alpha)-\underset{x^\alpha\in \Sm_{Z,m}}\sum wt(x^\alpha)+wt(x_l^m).
$$
By using \cite{HHL} again, we have
\begin{align*}
\mu([Y]_m^\star,\rho')&=\frac{mP_Y(m)}{n-l+1}\sum_{i=l}^{n}r_i-\sum_{x^\alpha\in \Sm_{Y,m}} wt(x^\alpha)\\
\mu([Z]_m^\star,\rho'')&=\frac{mP_Z(m)}{l+1}\sum_{i=0}^{l}r_i-\sum_{x^\alpha\in \Sm_{Z,m}} wt(x^\alpha)
\end{align*}
and plugging these in the Equation~(\ref{E:HHL}) produces the desired formula.
\end{proof}
In terms of monomial weights, Proposition~\ref{P:mu-reduction} takes the following form:
\begin{equation}\label{E:mu-reduction1}
\mu([X]_m^\star, \rho) = - \sum_{x^\a\in \Sm_{Y,m}} wt_{\rho'}(x^\a) - \sum_{x^\a\in \Sm_{Z,m}} wt_{\rho''}(x^\a) + \frac{mP(m)}{n+1} + mr_l.
\end{equation}
Using the lemma inductively as in the proof of Corollary~\ref{C:in}, we obtain the general case Proposition~\ref{P:mu-reduction} immediately.

\section{Decomposition of the barycenter}
Let $H$ be a hyperplane in $\bR^{n+1}$ defined by the equation $\sum_{i=0}^{n}x_i=m$ for some $m\in\bR$. For any given sequence of integers $n_{0}=0 < n_1 < \cdots < n_\ell = n$ and a subset $\{m_1,m_2,\cdots,m_\ell\}$ of $\bR$ such that $m=m_1+m_2+\cdots+m_\ell$, we aim to show that  $H$  can be decomposed into a sum of affine subspaces $H_i\subset \bR^{n+1}$, $i=1, \dots, \ell$, defined by

\begin{flalign*}
H_1=& \left\{{\bf a} \in\bR^{n+1}\,\left|\, \sum_{j=0}^{n_1}a_j=m_1; a_j =0 , n_1+1 \le j \le n \right.\right\}\\
H_i=& \left\{{\bf a}\in\bR^{n+1}\,\left|\, \sum_{j=n_{i-1}}^{n_{i}}a_j=m_i; a_k = 0, 0 \le j \le n_{i-1}-1, n_i+1 \le j \le n \right.\right\} \\
&\textrm{for $i=2,3,\dots,\ell-1$}\\
H_\ell=& \left\{{\bf a}\in\bR^{n+1}\,\left|\, \sum_{j=n_{\ell-1}}^{n}a_j=m_\ell; a_j = 0, 0 \le j \le n_{\ell-1}-1 \right.\right\}
\end{flalign*}

\begin{proposition}\label{P:bc-decomposition}
Let $H$ and $H_i$ be given as above, then $H=\sum_{i=1}^{\ell}H_i$. Moreover, any point of $H$ has a unique decomposition as a sum of elements in $H_i$.
\end{proposition}
\begin{proof}
Since the sum of the coordinate of any point in $H_i$ is $m_i$, the Minkowski sum $\sum_{i=1}^{\ell}H_i$ is contained in $H$. Conversely, for any given point $p=(p_0,p_1,\cdots,p_n)$, we can decompose $p$ as follows.
First, let
\[
q_{1j}=\left\{
\begin{array}{ll}
p_j & \textrm{if $0\leq j\leq n_{1}-1$}\\
m_1-\sum_{i=0}^{n_{1}-1}p_{i} & \textrm{if $j=n_1$}\\
0 & \textrm{else}\\
\end{array} \right.
\]
and note that $q_1=(q_{10},q_{11},\cdots,q_{1n})$ is an element of $H_1$.  Now, consider $p-q_1$ as an element of the hyperplane in $\{x \in \bR^n \, | \, x_0 = x_1 = \cdots = x_{n_1-1} = 0\} \simeq \bR^{n-n_1+1}$ defined by $\sum_{i=n_1}^{n}x_i=m-m_1$. In the same manner, we find $q_2\in H_2$ such that $p-q_1=q_2+(p-q_1-q_2)$ such that $p - q_1 - q_2$ is in the hyperplane of $\{x_0=x_1=\cdots = x_{n_2-1} =0\}$ defined by  $\sum_{i=n_2}^{n}x_i=m-m_1-m_2$. It is plain that repeating this procedure inductively produces a decomposition $p=q_1+q_2+\cdots+q_\ell$ where $q_i\in H_i$.

It remains to show the uniqueness of the decomposition. Let $p$ be a point in $H$ and suppose that $p=q_1+q_2+\cdots+q_\ell=r_1+r_2+\cdots+r_\ell$ where $q_i,r_i\in H_i$ for all $i=1,2,\dots,\ell$. Rearranging the terms, we have
\[
q_1-r_1=(r_2+r_3+\cdots+r_n)-(q_2+q_3+\cdots+q_n)
\]
which implies that $(q_1)_i=(r_1)_i$ for $0\leq i\leq n_{1}-1$  because the $i$th coordinate of the right hand side is zero when $i<n_1$. But since $\sum_{i=0}^{n}(q_1)_i=\sum_{i=0}^{n}(r_1)_i=m_1$ and $(q_1)_i = (r_1)_i = 0$ for $i > n_1$, it follows that $(q_1)_{n_1} = (r_1)_{n_1}$ and hence $q_1 = r_1$. In the same manner, one can easily show that $q_i=r_i$ for all $i\ge 2$ and this completes the proof.
\end{proof}

\begin{corollary}\label{C:bc-decomposition}
Retain notations $H_i$ from Proposition~\ref{P:bc-decomposition}. Let $P_i$ be a polytope contained in $H_i$ for $i=1,2,\dots,\ell$. We say that $H_i$ is the \emph{supporting hyperplane} for $P_i$. Let $p \in H$ and  $p=q_1+q_2+\cdots+q_\ell$ be the decomposition of $p$ as a sum of elements in $H_i$ given by Proposition~\ref{P:bc-decomposition}.
Then any $p$ is contained in $\sum_{i=1}^{\ell}P_i$ if and only if each $q_i$ is contained in $P_i$.
\end{corollary}
\begin{proof}
Suppose that $p$ is contained in $\sum_{i=1}^{\ell}P_i$. Then by definition of the Minkowski sum, there are elements $r_i\in P_i\subset H_i$ such that $p=r_1+r_2+\cdots+r_\ell$. But by the uniqueness of the decomposition, $r_i=q_i$ for all $i=1,2,\dots,\ell$. The other implication holds trivially.
\end{proof}

In the case of the state polytopes, as in Theorem~\ref{T:main}, the state polytopes
$$P_i := \st_m(I_{X_i}\cap k[x_{n_{i-1}},\dots, x_{n_i}])$$
are supported by the hyperplanes
$$H_i := \left\{x_0 = \cdots = x_{n_{i-1}-1} = 0, x_{n_i+1} = \cdots = x_n = 0, \sum_{i=0}^n x_i = m Q_i(m)\right\}$$
where $Q_i(m) = \dim_k (I_{X_i})_m$. Likewise, the state polytope of $X$ is supported by $H = \{ \sum_{i=0}^n x_i = m Q(m)\}$, $Q(m) = \dim_k (I_X)_m$. Let $\gamma$ be the barycenter of the state polytope of $X$ and $\tau$ be as in Remark~\ref{rem:main}. Applying Corollary~\ref{C:bc-decomposition}, we obtain
\begin{corollary}
Let $\gamma-\tau=q_1+q_2+\dots+q_\ell$ be the decomposition into a sum of elements in the supporting hyperplanes $H_i$. Then $\gamma$ is contained in $\st_m(I_X)$ if and only if each $q_i$ is contained in $\st_m( I_{X_i}\cap k[x_{n_{i-1}},\cdots,x_{n_i}])$.
\end{corollary}

\section{GIT of HIlbert points of pluricanonical curves}\label{S:applications}
Our main application is to the study of GIT of pluricanonical curves.

\subsection{Bi-canonical elliptic bridge} We revisit the state polytope analysis in \cite[Example~8.4]{MS}. Morrison and Swinarski considers the state polytope of a genus five curve of the form $C = \mathcal W_2 \cup E \cup \mathcal W_2$ where $\mathcal W_g$ denotes the Wiman curve of genus $g$ (Section~6.2, ibid) and $E$ is the elliptic curve $y^2 = x^3 - x$. According to the direct computation using the ideal of $C$,  it has $500,094$ initial ideals. By using the decomposition formula and Macaulay 2, we can compute its state polytope rather easily since the state polytopes of $\mathcal W_2$ and $E$ are fairly small.
Here, we give a complete description of the second state polytope.

\[
\txt{$\st_2(I_C) = \tau + $ }
\txt{$\quad \quad \quad \textup{Conv} \quad \quad \quad$\\
 $\overbrace{(2, 2, 1, 1, 2, {\bf 0}_7)}$ \\
 $(2, 2, 0, 3, 1, {\bf 0}_7)$ \\ $(2, 1, 1, 4, 0, {\bf 0}_7)$ \\ $(2, 0, 3, 3, 0, {\bf 0}_7)$ \\ $(0, 2, 3, 3, 0, {\bf 0}_7)$ \\ $(0, 0, 3, 4, 1, {\bf 0}_7)$ \\ $(0, 3, 1, 4, 0, {\bf 0}_7)$ \\ $(0, 1, 1, 5, 1, {\bf 0}_7)$ \\ $(0, 4, 0, 3, 1, {\bf 0}_7)$ \\ $(0, 2, 0, 4, 2, {\bf 0}_7)$ \\ $(2, 1, 3, 0, 2, {\bf 0}_7)$ \\ $(2, 0, 4, 1, 1, {\bf 0}_7)$ \\ $(0, 2, 4, 1, 1, {\bf 0}_7)$ \\ $(0, 0, 4, 2, 2, {\bf 0}_7)$ \\ $(0, 3, 3, 0, 2, {\bf 0}_7)$ \\ $(0, 1, 3, 1, 3, {\bf 0}_7)$ \\ $(0, 4, 1, 1, 2, {\bf 0}_7)$ \\
$\underbrace{(0, 2, 1, 2, 3, {\bf 0}_7)}$\\
}
\txt{$+$}
 \txt{$\quad \quad \quad \textup{Conv} \quad \quad \quad$\\
 $\overbrace{({\bf 0}_4, 2, 0, 1, 1, {\bf 0}_4)}$ \\
$({\bf 0}_4, 1, 0, 3, 0, {\bf 0}_4)$ \\
$({\bf 0}_4, 0, 0, 3, 1, {\bf 0}_4)$ \\
$({\bf 0}_4, 1, 0, 1, 2, {\bf 0}_4)$ \\
$({\bf 0}_4, 1, 2, 0, 1, {\bf 0}_4)$ \\
$\underbrace{({\bf 0}_4, 0, 2, 2, 0, {\bf 0}_4)}$\\
}
\txt{$+$}
\txt{$\quad \quad \quad \textup{Conv} \quad \quad \quad$\\
 $\overbrace{({\bf 0}_7, 2,1,1,2,2)}$ \\
$({\bf 0}_7, 3, 2, 1, 2, 0)$ \\ $({\bf 0}_7, 2, 1, 1, 4, 0)$ \\ $({\bf 0}_7, 1, 3, 0, 4, 0)$ \\ $({\bf 0}_7, 2, 4, 0, 2, 0)$ \\ $({\bf 0}_7, 1,  5, 1, 1, 0)$ \\ $({\bf 0}_7, 0, 4, 1, 3, 0)$ \\  $({\bf 0}_7, 1, 3, 0, 2, 2)$ \\ $({\bf 0}_7, 0, 4, 1, 1, 2)$ \\ $({\bf 0}_7, 2, 0, 3, 3, 0)$ \\ $({\bf 0}_7, 1, 1, 4, 2, 0)$ \\ $({\bf 0}_7, 3, 1, 3, 1, 0)$ \\ $({\bf 0}_7, 2, 0, 3, 1,  2)$ \\ $({\bf 0}_7, 1, 1, 4, 0, 2)$ \\ $({\bf 0}_7, 2, 2, 4, 0, 0)$ \\ $({\bf 0}_7, 1, 4, 3, 0, 0)$ \\ $({\bf 0}_7, 0, 3, 3, 2, 0)$ \\  $\underbrace{({\bf 0}_7, 0, 3, 3, 0, 2)}$
}
\]
The 2nd state polytope of $C$ is the $\tau=(7,7,7,7,4,8,8,4,7,7,7,7)$-translate of Minkowski sum as given above, where ${\bf 0}_i = (\underbrace{0, \dots, 0}_i)$.  The three columns are the (2nd) state polytopes of the three components $\mathcal W_2$, $E$ and $\mathcal W_2$.
Moreover, by Corollary~\ref{C:vertex}, the vertices of the Minkowski sum are precisely the sums of three vertices obtained by choosing one from each column. Hence the 2nd state polytope of $C$ has $18\cdot 18\cdot 6 = 1944$ vertices.

 Using Proposition~\ref{P:bc-decomposition}, we consider the decomposition of the barycenter. The barycenter of $\st_2(I_C)$ is $\gamma=(\frac{25}{3},\frac{25}{3},\dots,\frac{25}{3})$ and  $\gamma-\tau$ can be uniquely decomposed into a sum of elements in the supporting hyperplanes as follows.
\[
\gamma-\tau=\left(\frac{4}{3},\frac{4}{3},\frac{4}{3},\frac{4}{3},\frac{8}{3},{\bf 0}_7\right)+\left({\bf 0}_4,\frac{5}{3},\frac{1}{3},\frac{1}{3},\frac{5}{3},{\bf 0}_4\right)+\left({\bf 0}_7,\frac{8}{3},\frac{4}{3},\frac{4}{3},\frac{4}{3},\frac{4}{3}\right).
\]
Using Macaulay~2, we verify that the summands are not contained in the 2nd state polytope of $W_2$ and $E$. Therefore, the bicanonical image of $C$ is $2$-Hilbert unstable.

\subsection{$4$-canonical curve with a cuspidal tail} Let $Y = D\cup_p R$ be a pseudo-stable curve \cite{Schubert} consisting of a nonsingular curve of genus $g-1$, $g \ge 2$, and a rational curve $R$ with a cusp meeting in a node.
 We revisit the computation in \cite[Lemma~3]{HM} where the $m$-Hilbert instability for any $m\ge 2$ of the $4$-canonical image of $Y$
\[
Y  \hookrightarrow \bP(\Gamma(\omega_Y^{\otimes 4})) \simeq \bP^n
\]
 is proved.
Here,  $n = 7(g-1)$ and the isomorphism is given by choosing sections (homogeneous coordinates) such that
\[
\begin{array}{l}
R \subset \{x_0 = \dots = x_{7g-11} = 0\} \\
D \subset \{x_{7g-9} = \dots = x_n = 0\}.
\end{array}
\]
Let $l = 7g-10$ and apply Proposition~\ref{P:mu-reduction}. Let $\rho$ be the one-parameter subgroup with weight $(4, 4, \dots, 4, 3, 2, 0)$. The total $\rho$-weight of the degree two monomials in $x_l, \dots, x_n$ that are not in $\ii_{\prec_\rho}(R)$ is $35$. The degree two monomials in $x_0, \dots, x_l$ not in $\ii_{\prec_\rho}(D)$ contribute $15g-22$ to the total weight. Lastly, $mP(m)/(n+1) = m(8m-1)(4g-5)$. Putting all these together through Proposition~\ref{P:mu-reduction} or the monomial weight version (\ref{E:mu-reduction1}),
we get
\[
\mu([Y]_2, \rho) = -(35 + 2\cdot 4\cdot (15g-22)) + 30\cdot(4g-5) + 2\cdot 4 = -1.
\]
Likewise,
\[
\mu([Y]_3, \rho) = -(77 + 3\cdot 4\cdot (23g-34)) + 3\cdot23\cdot(4g-5) + 3\cdot 4 = -2
\]
from which it is deduced that $\mu([Y]_m,\rho) = -m+1$.

\subsection{Open rosaries}\label{S:rosary}
 We revisit the Hilbert-Mumford index computation of the special curves called `{\it rosaries}'. Recall from \cite[Definition~6.1]{HH2} that an open rosary of genus $r$ is $R = L_1\cup_{a_1} L_2\cup_{a_2} \cdots \cup_{a_r} L_{r+1}$ where $L_i$ are smooth genus zero curves and $a_i$'s are tacnodes. Each $L_i$ intersects with $L_j$ if and only if $|i-j|=1$.  In \cite[Section~8.1]{HH2}, the Hibert-Mumford index of an open rosary (with respect to a 1-ps coming from its automorphism group)  embedded by $\omega_R^{\otimes 2}(2a_0+2a_{r+1})$ is computed, where $a_0$ and $a_{r+1}$ are smooth points of $L_1$ and $L_{r+1}$, respectively.
\begin{figure}[!htb]\labellist \small\hair 2pt
\pinlabel $a_0$ at 215 594
\pinlabel $a_1$ at 252 616
\pinlabel $a_2$ at 290 616
\pinlabel $a_3$ at 325 616
\pinlabel $a_4$ at 346 594
\endlabellist
  \begin{center}
    \includegraphics[width=2.5in]{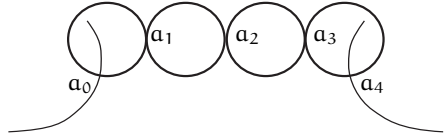}
  \end{center}
   \caption{An open rosary of genus three} \label{F:open-rosary4}
\end{figure}
 Note that since $a_1, \dots, a_r$ are tacnodes, we are not able to apply the decomposition formula directly. However, in the proposition below, we shall demonstrate that a similar argument can be used to obtain a systematic analysis of the initial ideal.


\begin{proposition} The initial ideal of $R$ with respect to the $\rho$-weighted Lex order satisfies the following decomposition: Let $T_l^d$ denote the set of degree $d$ monomials in $x_0, \dots, x_{3l+2}$ (resp. $x_0, \dots, x_{3r}$) which involve $x_i$ and $x_j$ for $i<3l-2$ and $j>3l-1$ for $l = 1, \dots, r-1$ (resp. $l = r$). We define $x_l = x_0$ if $l < 0$ and $x_l = x_{3r}$ if $l > 3r$.
\begin{enumerate}
\item[(i)] \textup{(degree 2 piece)}
\[
(\ii I_R)_2 \cup \{x_{3l-2}^2\}_{l=1}^{r}  = \cup_{l=1}^{r+1} \ii \left(I_{L_l}\cap k[x_{3l-5}, \dots, x_{3l-1}]\right)_2 \cup (\cup_{l=1}^r T_l^2)
\]
\item[(ii)] \textup{(degree 3 piece)}
\[
(\ii I_R)_3 \cup \{x_{3l-2}^3,x_{3l-2}^2x_{3l-1}\}_{l=1}^{r}  =  \cup_{l=1}^{r+1}\ii  \left(I_{L_l}\cap k[x_{3l-5}, \dots, x_{3l-1}]\right)_3 \cup (\cup_{l=1}^r T_l^3)
\]
\end{enumerate}
\end{proposition}

\begin{proof} The proof follows the idea of Lemma~\ref{L:in} closely, but more care needs to be exercised because of the additional overlapping coordinates.

\noindent (i) \, (degree 2 piece) \, Suppose that $x^{\a}\in(\ii I_R)_2$, then $x^{\a}=\ii f$ for some homogeneous quadratic $f\in I_R$. If $x^{\a}$ is not contained in $\cup_{l=1}^r T_l$, then $x^{\a}$ is contained in $k[x_{3l-5},\dots,x_{3l-1}]$ for some $l$. Let $g=f(0,\dots,0,x_{3l-5},\dots,x_{3l-1},0,\dots,0)$. Then $g\in I_{L_l}\cap k[x_{3l-5},\dots,x_{3l-1}]$ and therefore, $x^{\a}$ is contained in $(\ii I_{L_l}\cap k[x_{3l-5},\dots,x_{3l-1}])_2$. Since $x_{3l-2}^2-x_{3l-1}x_{3l}\in I_{L_{l+1}}$, we have $x_{3l-2}^2\in\ii I_{L_l+1}$.

For the other inclusion, suppose that $x^{\a}\in(\ii I_{L_l}\cap k[x_{3l-5},\dots,x_{3l-1}])_2$ for some $l=1,2,\dots,r+1$. If $l=1$ or $r+1$, then $x^{\a}=x_0x_2$ or $x_{3r-2}^2$ respectively, but we know that $x_0x_2=\ii(x_0x_2-x_1^2-x_2x_3)\in\ii I_{R}$. Hence $x^{\a}$ is contained in the left hand side. If $l=2,3,\dots,r$, then
\[
\begin{array}{l}
I_{L_l}=\langle x_{3l-2}^2-x_{3l-3}x_{3l-1},x_{3l-3}x_{3l-2}-x_{3l-5}x_{3l-1},x_{3l-5}x_{3l-2}-x_{3l-4}x_{3l-1},\\
x_{3l-3}^2-x_{3l-4}x_{3l-1},x_{3l-5}x_{3l-3}-x_{3l-4}x_{3l-2},x_{3l-5}^2-x_{3l-4}x_{3l-3}\rangle.
\end{array}
\]
Let's denote the generators of $I_{L_l}$ by $f_i$ for $i=1,2,\dots,6$ in the order shown above.\\
The initial ideal of $L_{l}$ is
\[
\begin{array}{lll}
\ii I_{L_l} & = & \langle x_{3l-3}x_{3l-1},x_{3l-4}x_{3l-1},x_{3l-4}x_{3l-2}^2,x_{3l-5}x_{3l-1},x_{3l-5}x_{3l-2},\\
 & & x_{3l-5}x_{3l-3},x_{3l-5}^2\rangle
\end{array}
\]
 This shows that $x^{\a}$ comes from the initial term of some generator of $I_{L_{l}}$. That is, $x^{\a}=\ii(f_i)$ for some $i$. But $f_i\in I_R$ for $i=2,3,4,5$, it suffices to consider the case $i=1,6$. For the case $i=1$, then $x^{\a}=x_{3l-3}x_{3l-1}$ and so $x^{\a}\in\ii(x_{3l-2}^2-x_{3l-3}x_{3l-1}-x_{3l-1}x_{3l})\in\ii I_{R}$. For $i=6$, then $x^{\a}=x_{3l-5}^2$. Obviously, $\cup_{l=1}^r T_l$ is contained in $\ii I_{R}$.

\

\noindent (ii) \, (degree 3 piece) \, Since $x_{3l-2}^2\in\ii I_{L_{l+1}}$, $x_{3l-2}^3$ and $x_{3l-2}^2x_{3l-1}$ are also elements in $I_{L_{l+1}}$. If $x^{\a}\in\ii I_R$ and not in $\cup_{l=1}^r T_l$, we can easily show that $x^{\a}$ is contained in the right hand side by the same proof of the degree 2 case.

For the other inclusion, suppose that $x^{\a}$ is a degree three monomial in $\ii(I_{L_l}\cap k[x_{3l-5},\dots,x_{3l-1}])$ for some $l$. If $l=1$, then $x^{\a}=x_i(x_0x_2)$ for some $i=0,1,2$. Hence $x^{\a}=\ii(x_i(x_0x_2-x_1^2-x_2x_3))\in\ii I_R$. If $l=r+1$, then $x^{\a}=x_i(x_{3r-2})^2$ for some $i=3r-2,3r-1,3r$. For $i=3r$, we have $x^{\a}=\ii(x_{3r}(x_{3r-2}^2-x_{3r-1}x_{3r}))\in\ii I_{R}$.

Suppose that $l=2,3,\dots,r$, then by using the same notation in degree 2 case, we know that $x^{\a}=\ii(x_if_j)$ for $j=1,2,\dots,6$ and $i=3l-5,\dots,3l-1$ except for the case $x^{\a}=x_{3l-4}x_{3l-2}^2$. If $j=2,3,4,5$, then $x_if_j\in I_{R}$ and so $x^{\a}\in\ii I_{R}$. If $j=1$, then $x^{\a}=\ii(x_i(x_{3l-2}^2-x_{3l-3}x_{3l-1}-x_{3l-1}x_{3l}))\in\ii I_{R}$. Hence it remains to consider the case $j=6$ and $x^{\a}=x_{3l-4}x_{3l-2}^2$.

For the case $j=6$, suppose that $x^{\a}=\ii(x_if_6)=x_ix_{3l-5}^2$. For $i=3l-5, 3l-4$, there is nothing to prove. If $i>3l-4$, then $x^{\a}=\ii(x_{3l-5}f_j)$ for $j=2,3,5$. Hence $x^{\a}\in\ii I_{R}$ since $f_j\in I_R$ where $j=2,3,5$.
For the case $x^{\a}=x_{3l-4}x_{3l-2}^2$, $x^{\a}=x_{3l-4}x_{3l-2}^2=\ii\{x_{3l-4}(f_1-x_{3l-1}x_{3l})+x_{3l-3}(f_4)\}\in\ii I_R$ since $f_1-x_{3l-1}x_{3l}$ and $f_4$ are elements in $I_R$. Finally, $\cup_{l=1}^r T_l$ is contained in $\ii I_{R}$. This completes the proof.
\end{proof}

By using the above proposition, we can compute inductively the sum of $\rho$-weights of degree 2 and 3 monomials in the initial ideal of $I_R$. Let $w_i(r)$ be the sum $\underset{x^{\a}\in\ii(I_{R})_i}\Sm wt_{\rho}(x^{\a})$ for $i=2,3$ where $R$ is the open rosary of length $r+1$. Then
\[
w_2(r)=\left\{
\begin{array}{ll}
w_2(r-2)+(72r-110) & \textrm{for odd $r>1$}\\
w_2(r-2)+(72r-92) & \textrm{for even $r>2$}\\
\end{array} \right.
\]
For the degree 3 case,
\[
w_3(r)=\left\{
\begin{array}{ll}
w_3(r-2)+(162r^2-162r-57) & \textrm{for odd $r>1$}\\
w_3(r-2)+(162r^2-108r-66) & \textrm{for even $r>2$}\\
\end{array} \right.
\]
Since $w_2(1)=6,w_2(2)=52,w_3(1)=34,w_3(2)=366$, we can easily compute the $w_i(r)$ as follows.
\[
w_2(r)=\left\{
\begin{array}{ll}
18r^2-19r+7 & \textrm{for odd $r$}\\
18r^2-10r & \textrm{for even $r$}\\
\end{array} \right.
\]
and
\[
w_3(r)=\left\{
\begin{array}{ll}
27r^3+\frac{81}{2}r^2-\frac{111}{2}r+22 & \textrm{for odd $r$}\\
27r^3+54r^2-33r & \textrm{for even $r$}\\
\end{array} \right.
\]
which is consistent with the results in \cite[Section~8.1]{HH2}.

\bibliographystyle{alpha}
\bibliography{st-decomposition}

\end{document}